\documentclass[12pt,dvips]{article}

\usepackage{amsmath,amsfonts,amssymb,amsthm,graphicx,verbatim,subfig}
\usepackage[all]{xy}

\ifx\pdftexversion\undefined

\usepackage[a4paper,colorlinks,link
=black,filecolor=black,citecolor=black,urlcolor=black,pdfstartview=FitH]{hyperref}
\else

\usepackage[a4paper,colorlinks,linkcolor=black,filecolor=black,citecolor=black,urlcolor=black,pdfstartview=FitH]{hyperref}
\fi

\font\sixbb=msbm6
\font\eightbb=msbm8
\font\twelvebb=msbm10 scaled 1095
\newfam\bbfam
\textfont\bbfam=\twelvebb \scriptfont\bbfam=\eightbb
                           \scriptscriptfont\bbfam=\sixbb

\newtheorem*{theorem*}{\bf Theorem}

\title{Chernoff's Inequality - A very elementary proof}
\begin{document}
\author{Nathan Linial\thanks{Department of Computer Science, Hebrew University, Jerusalem 91904,
    Israel. e-mail: nati@cs.huji.ac.il~. Supported by ISF and BSF grants.}
  \and {Zur Luria\thanks{Department of Computer Science, Hebrew University, Jerusalem 91904,
    Israel. e-mail: zluria@cs.huji.ac.il~.}} }

\date{}

\maketitle
\pagestyle{plain}

\begin{abstract}
We give a very simple proof of a strengthened version of Chernoff's Inequality. We derive the same conclusion from much weaker assumptions.
\end{abstract}

\section*{The theorem}

\begin{theorem*} \label{lp_theorem}
Let $X_1,\ldots,X_n$ be indicator random variables, $0 < \beta < 1$ and $ 0 < k < \beta n $.
Then
\begin{equation}\label{main}
\Pr \left( \sum_{i=1}^n{X_i}\geq \beta n\right) \leq
\frac{1}{\binom{\beta n}{k}}\sum_{|S|=k}{\Pr\left(\wedge_{i \in S} (X_i = 1) \right)}.
\end{equation}
In particular, if
$$
\Pr\left(\wedge_{i \in S} (X_i = 1) \right) \leq \alpha^k
$$
for every $S$ of size $k = \left(\frac{\beta-\alpha}{1-\alpha}\right)n$, where $0 < \alpha < \beta$ then
$$
\Pr \left( \sum_{i=1}^n{X_i}\geq \beta n\right)  \leq e^{-D(\beta||\alpha) n} .
$$
\end{theorem*}

If one makes the stronger assumption that $\Pr(\wedge_{i \in S} \left(X_i = 1\right)) \leq \alpha^{|S|}$ for {\em every} $S\subseteq\{1,\ldots,n\}$, this reduces to Theorem 1.1 from \cite{IK10}. Under the even more restrictive assumption that the $X_i$ are i.i.d., this reduces to the usual statement of Chernoff's Inequality.

\begin{proof}Let $w_S = \Pr(X_i = 1 \Leftrightarrow i \in S)$. Then
$$
\Pr(\wedge_{i \in S} \left(X_i = 1\right)) = \sum_{T: S \subseteq T}{w_T}
$$
and
$$
\Pr \left( \sum_{i=1}^n{X_i}\geq \beta n\right) = \sum_{T:|T| \geq \beta n}{w_T}.
$$
Since $\binom{x}{k}$ is an increasing function of $x$ and since $w_S \geq 0$, we have
\[
\sum_{|T|\ge\beta n}{w_T} \le
\frac{1}{\binom{\beta n}{k}}\sum_T {|T|\choose k}w_T =
\frac{1}{\binom{\beta n}{k}} \sum_{|S|=k}{ \sum_{T: S \subseteq T}{w_T} } =
\frac{1}{\binom{\beta n}{k}} \sum_{|S|=k}{\Pr(\wedge_{i \in S} \left(X_i = 1\right))},
\]
as claimed. For the second part of the theorem note that if we assume that
$\Pr\left(\wedge_{i \in S} (X_i = 1) \right) \leq \alpha^k$
for every $S$ of size $k = \left(\frac{\beta-\alpha}{1-\alpha}\right)n$ then the above inequality becomes
$$
\Pr \left( \sum_{i=1}^n{X_i}\geq \beta n\right) \leq
\frac{\binom{n}{k}}{\binom{\beta n}{k}}\alpha^k \leq e^{-D(\beta||\alpha)n}
$$

where the last inequality follows from standard entropy estimates of binomials.
\end{proof}

It is easy to see that inequality~\ref{main} is tight. It holds with equality iff $\Pr(\sum X_i=\beta n\text{~or~}\sum X_i < k)=1$.

\end{document}